\documentclass[12pt,a4paper]{article}
\usepackage[margin=1in]{geometry}  
\usepackage{setspace}
\usepackage{mathptmx} 
\usepackage{amsmath,amssymb,amsthm,mathtools}
\usepackage{bm}
\usepackage{graphicx}
\usepackage{booktabs}
\usepackage{enumitem}
\usepackage{hyperref}
\usepackage[nameinlink,noabbrev]{cleveref}
\usepackage{microtype}
\usepackage{algorithm}
\usepackage{algorithmic}
\usepackage{natbib}
\usepackage{multirow}

\theoremstyle{plain}
\newtheorem{theorem}{Theorem}
\newtheorem{lemma}[theorem]{Lemma}
\newtheorem{proposition}[theorem]{Proposition}
\newtheorem{corollary}[theorem]{Corollary}
\theoremstyle{definition}
\newtheorem{definition}[theorem]{Definition}
\newtheorem{assumption}{Assumption}
\theoremstyle{remark}
\newtheorem{remark}[theorem]{Remark}
\newtheorem{example}[theorem]{Example}
\newcommand{\E}{\mathbb{E}}
\newcommand{\R}{\mathbb{R}}
\newcommand{\KL}{\mathrm{KL}}
\newcommand{\MI}{I}
\newcommand{\Var}{\operatorname{Var}}

\begin{document}

\title{Entropy-Regularized Inference:\\
A Predictive Approach}

\author{%
Nicholas G.\ Polson\\
{\normalsize Booth School of Business, University of Chicago}\\
{\normalsize \texttt{ngp@chicagobooth.edu}}
\and
Daniel Zantedeschi\\
{\normalsize Muma College of Business, University of South Florida}\\
{\normalsize \texttt{danielz@usf.edu}}%
}

\date{}

\maketitle

\begin{abstract}
\noindent
Predictive inference requires balancing statistical accuracy against informational
complexity, yet the choice of complexity measure is usually imposed rather than derived. We treat econometric objects
as \emph{predictive rules}, mappings from information to reported predictive
distributions, and impose three structural requirements on evaluation: locality,
strict propriety, and coherence under aggregation (coarsening/refinement) of outcome
categories. These axioms characterize (uniquely, up to affine transformations) the logarithmic score and induce Shannon mutual information (Kullback--Leibler divergence) as the
corresponding measure of predictive complexity. The resulting entropy-regularized
prediction problem admits Gibbs-form optimal rules, and we establish an essentially
complete-class result for the admissible rules we study under joint risk--complexity dominance. Rational
inattention emerges as the constrained dual, corresponding to frontier points with
binding information capacity. The entropy penalty contributes additive curvature
to the predictive criterion; in weakly identified settings, such as weak instruments
in IV regression, where the unregularized objective is flat, this curvature
stabilizes the predictive criterion. We derive a local
quadratic (LAQ) expansion connecting entropy regularization to classical
weak-identification diagnostics.


\end{abstract}

\onehalfspacing

\section{Introduction}
\label{sec:intro}

Predictive inference requires controlling informational complexity. When a predictive
rule refines its reports too sharply in response to states or signals, the resulting
predictions become unstable, sensitive to noise, overfitted to training data, or
ill-posed when identification is weak. The natural question is: what penalty on
predictive refinement is coherent with how predictions are evaluated? This paper
shows that once evaluation satisfies minimal coherence requirements, a specific
penalty is characterized within our framework: entropy (equivalently, Kullback--Leibler divergence).
Unlike Bayesian shrinkage, entropy regularization operates in predictive space
and does not rely on prior mass to repair flat likelihood directions.
The resulting regularized prediction rules take a Gibbs (exponential tilting) form,
and these rules form an essentially complete class within the admissible set
under a risk--complexity ordering.

The mechanism is geometric. Entropy regularization penalizes the mutual information
between states and reports, which measures how much the predictive distribution
varies with the information available. This penalty contributes curvature to the
predictive criterion: the regularized objective has Hessian
$H_\lambda(\theta) = H_U(\theta) + \lambda^{-1} H_{\MI}(\theta)$, where $H_U$ is the
utility curvature and $H_{\MI}$ is the information curvature. The key property is
\emph{additivity}: because the logarithmic score decomposes additively under
coarsening and refinement of outcome categories, information is composable across
partitions. This additivity implies exponential tilting for the optimal rule under the induced objective (the Gibbs form) as the
structure of optimal predictive rules, not as an assumption but as a consequence
of coherent evaluation. This becomes operationally decisive precisely when the
unregularized criterion is nearly flat.

Weak identification provides a natural econometric setting where this principle
matters. When the criterion linking parameters to observables lacks curvature, as
in weak-instrument IV, sparse discrete choice, or flat likelihood surfaces, small
data perturbations induce large changes in fitted predictions
\citep{staiger1997instrumental,stock2002survey,andrews2012estimation}. The entropy
penalty adds curvature in these directions: $\lambda^{-1} H_{\MI}(\theta)$ has mass in
\emph{all} directions, including those where $H_U(\theta)$ is deficient. The
result is a well-posed predictive map even when the unregularized criterion is
nearly singular. From this viewpoint, weak identification is not the premise but
the consequence: it is the regime where coherent regularization of predictive
refinement becomes essential.

\paragraph{Contributions}
The paper delivers four main results:
\begin{enumerate}[label=(\roman*),nosep]
\item \emph{Uniqueness of the log score and Gibbs form.} We impose three axioms
on predictive evaluation (locality, strict propriety, and aggregation coherence) and
show these uniquely select the logarithmic score. The induced complexity measure is
KL divergence, and optimal rules take the Gibbs (exponential tilting) form
(Sections~\ref{sec:setup}--\ref{sec:gibbs}).

\item \emph{LAQ interpretation and curvature restoration.} We derive a local
asymptotic quadratic (LAQ) expansion showing how the entropy penalty contributes
curvature in directions where the utility Hessian is flat. This connects entropy
regularization to classical weak-identification diagnostics and explains why
prediction stabilizes under weak identification
(Section~\ref{sec:identification}).

\item \emph{Complete-class characterization.} Under a joint ordering preferring
higher expected utility and lower information complexity, Gibbs rules form an
essentially complete class within the aggregation-coherent, locally scored
predictive rules we study (in the sense of weak convergence on finite outcome
spaces). Every non-Gibbs rule is dominated by some Gibbs rule achieving equal
or better predictions with equal or lower complexity. Rational inattention
emerges as the constrained dual (Section~\ref{sec:complete}).

\item \emph{Connections to entropy econometrics.} We position the framework
relative to empirical likelihood, GEL/ETEL, Bayesian empirical likelihood,
rational inattention, and predictive density aggregation
(Section~\ref{subsec:literature}).
\end{enumerate}

\paragraph{Roadmap}
Section~\ref{subsec:literature} reviews related literature.
Sections~\ref{sec:setup}--\ref{sec:coherence} define the setup and derive uniqueness of the log score.
Section~\ref{sec:gibbs} characterizes Gibbs optimal rules.
Section~\ref{sec:complete} establishes the complete-class theorem.
Section~\ref{sec:identification} develops the curvature decomposition and LAQ expansion.
Section~\ref{sec:numerical} provides Monte Carlo evidence.
Section~\ref{sec:discussion} concludes.
Full proofs appear in \ref{app:proofs} and \ref{app:complete}; additional simulations in \ref{app:simulations}.

\subsection{Related Literature and Positioning}
\label{subsec:literature}

This paper connects entropy and divergence methods in econometrics to proper
scoring rules, rational inattention, and weak identification by making an
\emph{econometric object choice} about what econometrics delivers. The object is
not an estimator or a parameter, but a \emph{predictive rule}: a mapping from
available information into a reported predictive distribution. From this
viewpoint, regularization is a restriction on the \emph{admissible modes of
predictive refinement}. The central question becomes: which penalties on
refinement are compatible with coherent predictive evaluation?

\medskip
\noindent\emph{Map of literatures and our contribution.}
We connect four literatures that usually interact only indirectly:
(i) entropy/divergence econometrics (EL/GEL/ETEL; BEL and moment-based Bayes),
(ii) proper scoring rules and forecast evaluation,
(iii) rational inattention and information design, and
(iv) weak identification and ill-conditioned criteria. Our contribution is not
a new divergence criterion or a new weak-IV confidence procedure. Rather, we
provide an axiomatic, prediction-first derivation of the KL/Shannon geometry as
the unique refinement penalty compatible with \emph{aggregation-coherent} forecast
evaluation, and we show how this geometry adds curvature in directions where the
unregularized predictive criterion degenerates. Our claims are conditional on the stated coherence axioms and the predictive-rule object; we do not argue that other objectives or penalty geometries are invalid, only that within this framework KL/entropy is uniquely characterized.

\subsubsection*{1. Divergence Methods in Econometrics: EL/GEL/ETEL and Bayesian Moment Criteria}

This paper connects directly to entropy and divergence methods in econometrics.
The conceptual entry point is exponential tilting: the Gibbs form
\eqref{eq:gibbs} is the predictive analogue of maximum entropy under moment
constraints. When an exchangeable sequence is conditioned on empirical moment
constraints, the resulting predictive law is an $I$-projection, an exponential
tilt of the reference measure, yielding an exponential-family structure
\citep{polsonzantedeschi2025definetti_sanov}. Relatedly,
\citet{polsonzantedeschi2025momentconditioning} develop a finite Bayesian
moment-conditioning framework that connects empirical likelihood, Bayesian
empirical likelihood, and GMM through a common predictive geometry. EL/GEL/ETEL and
related methods operationalize these divergence ideas in estimation.

\paragraph{Empirical likelihood, GEL, and ETEL}
Empirical likelihood \citep{owen2001empirical} and generalized empirical
likelihood (GEL) methods \citep{kitamura1997nonparametric,newey2004higherorder}
define estimators via divergence minimization subject to moment restrictions.
Exponentially tilted empirical likelihood (ETEL) \citep{schennach2007etel} uses
KL divergence specifically and exhibits favorable higher-order properties under
standard regularity; this is the frequentist side of KL tilting. On the Bayesian
side, parallel constructions replace the likelihood with empirical-likelihood
criteria built from moment restrictions, yielding Bayesian empirical likelihood
(BEL) and related pseudo-posterior approaches
\citep{lazar2003bayesianEL,chamberlain2003nonparametric,chibshinsimoni2018}. These methods typically
\emph{choose} a divergence for computational, robustness, or higher-order
reasons. ETEL, for instance, adopts KL for higher-order frequentist optimality. We show
that KL is uniquely selected once predictive evaluation satisfies aggregation
coherence, so the divergence plays fundamentally different roles: asymptotic
in ETEL, axiomatically derived here.

\paragraph{Information divergence as regularization}
Information-theoretic methods have a long history in econometrics,
including early perspectives on
Bayesian information processing \citep{zellner1988information}, syntheses on
entropy econometrics \citep{golan2002ieeeditor}, and applications linking entropy
to predictability in financial returns \citep{maasoumi2002entropy}. Divergence
penalties also appear in robust Bayesian analysis \citep{berger1994overview},
ambiguity-averse approaches to model uncertainty \citep{hansen2008fragile}, and
distributionally robust optimization. Our contribution clarifies the predictive
interpretation: entropy/KL regularization is the unique penalty compatible with
coherent predictive evaluation.

\subsubsection*{2. Proper Scoring Rules and Aggregation-Coherent Evaluation}

Proper scoring rules provide the decision-theoretic foundation for evaluating
probabilistic forecasts \citep{gneiting2007strictly,giacomini2006tests}, and the
logarithmic score is a workhorse for density forecast evaluation
\citep{diebold1998evaluating,amisano2007comparing,geweke2010comparing}. Local
scoring rules are natural in econometrics because they evaluate a report using
only the probability assigned to the realized outcome. A classical
characterization selects the log score as the unique strictly proper local score
(up to affine transformations) \citep{bernardo1979expected}, with modern
refinements of locality conditions \citep{dawid2012proper}.

Our contribution elevates a further coherence requirement,
\emph{aggregation coherence}, invariance under coarsening and refinement of
outcome partitions, to a structural axiom. This requirement is essential for
applied econometrics: outcomes are routinely grouped (income brackets), censored
(top-coded earnings), or discretized (choice categories), and evaluation must be
consistent across resolutions. Aggregation coherence rules out Bregman and
$f$-divergence alternatives that fail global coherence
\citep{BloedelDentiPomatto2025,FosgerauEtAl2020}, and uniquely selects Shannon
entropy and KL divergence as the geometry of predictive refinement. In particular,
this axiom explains why ``entropy regularization'' in our sense is not merely one
penalty among many: the KL/Shannon geometry is the one compatible with coherent
evaluation across resolutions.

\subsubsection*{3. Rational Inattention, Information Design, and the Dual Interpretation}

Rational inattention posits that decision makers face information capacity
constraints measured by Shannon mutual information
\citep{sims2003implications,matejka2015rational,mackowiak2009optimal}.
\citet{CaplinDeanLeahy2019} develop the connection between rational inattention,
optimal consideration sets, and stochastic choice, showing how information costs
shape discrete decisions. Related work emphasizes the choice of information
structures and their welfare implications \citep{kamenica2011bayesian,caplin2015revealed}.
In most applications, the Shannon cost is postulated. Our paper provides a
coherence-based derivation: mutual information arises because it is the
refinement measure uniquely induced by coherent local scoring. Rational
inattention appears as the constrained dual of entropy-regularized prediction
(see Section~\ref{sec:complete}).

\paragraph{Predictive density aggregation}
Recent work studies regularized mixtures for aggregating predictive densities,
with applications to macroeconomic forecasting \citep{diebold2023aggregation}.
Related work develops proper scoring approaches for interval and quantile
forecasts \citep{askanazi2018interval}. For empirical comparisons between
aggregation strategies in ill-conditioned regression and entropy-based methods,
see \citet{tavares2025aggregation}.

\subsubsection*{4. Weak Identification: Curvature Degeneracy and Predictive Stabilization}

Weak identification occurs when the criterion linking parameters to observables
provides insufficient curvature to localize the parameter, producing nonstandard
asymptotics and unreliable Wald-type inference
\citep{staiger1997instrumental,stock2002survey,andrews2012estimation,andrews2019weak}.
A substantial literature develops identification-robust procedures and confidence
sets that remain valid when curvature degenerates
\citep{moreira2003clr,mikusheva2010robust,kleibergen2005testing}. Two recent lines
of work are especially relevant. On the geometric side,
\citet{andrewsmikusheva2016,andrewsmikusheva2022} characterize weak identification
through tangent cones and conditional inference, delivering LAMN limit experiments
and optimal decision rules. On the GMM side,
\citet{antoine2009efficient,antoine2012efficient,antoine2020testing} develop
asymptotics under nearly-weak identification and propose identification-strength
tests for intermediate regimes. In Section~\ref{sec:identification} we show that
our curvature decomposition connects these programs; we identify the bridge but
defer the full distributional theory to future work.

Our contribution shows that entropy regularization, derived from coherence, stabilizes the predictive criterion under curvature degeneracy: the entropy term contributes additive curvature, yielding a well-conditioned local quadratic approximation in directions where the unregularized utility curvature is deficient. We do not propose a new identification-robust confidence procedure; we show, within our framework, why the regularized predictive criterion remains well-conditioned when the unregularized criterion flattens, and we provide an LAQ diagnostic bridge. We focus on prediction and criterion geometry; developing full frequentist optimal inference procedures is beyond scope.

\paragraph{Running example: Weak instruments.}
We illustrate the theory throughout with instrumental variables regression. In
the standard setup $Y_i = X_i \beta + \varepsilon_i$, $X_i = Z_i'\Pi + V_i$,
weak instruments ($\Pi \approx 0$) create flat directions in the predictive
criterion. The entropy penalty adds curvature precisely there, softening the
predictive distribution $p(Y \mid X, Z)$ rather than shrinking coefficients
toward zero. The object of interest is the prediction, not the parameter.

\paragraph{Why standard Bayesian approaches do not resolve weak IV}
Finding valid instruments is hard; \citet{rossi2014ivmarketing} argues that
instruments used in applied work are rarely valid and that IV methods exhibit
poor finite-sample properties including substantial bias and high sampling
variance. Standard Bayesian approaches do not automatically resolve weak
identification. \citet{lopes2014bayesianiv} show that Jeffreys priors yield
pathological posteriors under weak instruments: the posterior concentrates on
implausible regions or exhibits extreme sensitivity to minor specification
changes. Similarly, \citet{kleibergen1998bayesianiv} and
\citet{hoogerheide2007conjugateiv} document that conjugate priors fail to
regularize the flat directions in the IV likelihood; the Jeffreys prior in
particular produces posteriors identical in functional form to the LIML sampling
distribution, inheriting its instability. \citet{sims2007thinkingiv} emphasizes
that prior beliefs play an ``inevitably strong role'' when instruments are weak,
and that decision-theoretic reasoning, rather than reference priors, is required
to obtain well-behaved inference. Entropy regularization provides precisely this:
a coherence-based penalty that stabilizes prediction without
relying on arbitrary prior specification.

\subsubsection*{5. Regularization in Econometrics and Statistical Learning: Why Penalty Geometry Matters}

Classical regularization stabilizes ill-conditioned objectives by adding
curvature: ridge regression \citep{hoerl1970ridge} penalizes $\|\beta\|_2^2$,
lasso \citep{tibshirani1996regression} penalizes $\|\beta\|_1$, and elastic net
blends the two \citep{zou2005regularization}. Bayesian shrinkage achieves
analogous stabilization through priors \citep{park2008bayesian,carvalho2010horseshoe},
with recent work on global-local shrinkage providing adaptive regularization
\citep{bhadra2016globallocal,polson2010shrink}. In econometrics, regularization is
central in ill-posed inverse problems and high-dimensional moment models
\citep{carrasco2012regularization,florenssimoni2012,hansen2016efficient}.

These contributions explain how penalties work and what guarantees they deliver,
but typically treat the geometry of the penalty as exogenous. Our point of
departure is that for prediction, penalty geometry can be derived from coherence
requirements on evaluation, yielding an information-theoretic (KL/entropy) notion
of complexity rather than a parameter-norm notion.

\section{Predictive Rules, Scoring, and Information Complexity}
\label{sec:setup}

The econometric object throughout is the \emph{predictive rule}, not an estimator or
parameter. The goal is to characterize how predictive distributions should respond
to information states when evaluated by coherent scoring rules.

Let $(\mathcal{X},\mathcal{F}_X)$ denote the \emph{outcome} space, equipped with a
$\sigma$-finite reference measure $\mu$ when densities are invoked. Let
$(\mathcal{S},\mathcal{F}_S)$ denote the \emph{information state} space (signals,
covariates, instruments, moments, or other summaries available to the analyst).
We write $P_S$ for a prior (or population) distribution on $(\mathcal{S},\mathcal{F}_S)$.

A \emph{predictive report} is a probability law on $(\mathcal{X},\mathcal{F}_X)$.
We let $\mathcal{Y}\subseteq \Delta(\mathcal{X})$ denote the set of admissible reports.
When $y\in\mathcal{Y}$ is absolutely continuous with respect to $\mu$, we write its
density as $y(x)=\frac{dy}{d\mu}(x)$.

A \emph{stochastic predictive rule} (a \emph{channel}) is a Markov kernel
\begin{equation}
\label{eq:kernel}
f:\mathcal{S}\to \Delta(\mathcal{Y}),\qquad s\mapsto f(\cdot\mid s),
\end{equation}
mapping each information state $s$ into a distribution over predictive reports. Given
$(P_S,f)$, the induced joint law on $(S,Y)$ is
\begin{equation}
\label{eq:joint}
P_{SY}(ds,dy)=P_S(ds)\,f(dy\mid s),
\end{equation}
with marginal distribution of reports
\begin{equation}
\label{eq:marginalq}
q(dy)=\int f(dy\mid s)\,P_S(ds).
\end{equation}

\begin{definition}[Scoring rule]
\label{def:scoring}
A scoring rule is a map $S:\mathcal{X}\times\mathcal{Y}\to\R$ that assigns a numerical
score $S(x,y)$ when outcome $x$ is realized and report $y$ was issued. The scoring
rule is \emph{proper} if for every distribution $P\in\Delta(\mathcal{X})$,
$\E_{X\sim P}[S(X,P)]\ge \E_{X\sim P}[S(X,Q)]$ for all $Q\in\mathcal{Y}$, and
\emph{strictly proper} if equality holds only when $Q=P$ (up to $P$-null sets).
\end{definition}

\begin{definition}[Locality]
\label{def:locality}
A scoring rule is \emph{local} if there exists a function $s:\R_+\to\R$ such that
whenever the report $y$ admits a density with respect to $\mu$,
$S(x,y)= s\!\left(\frac{dy}{d\mu}(x)\right)$.
Thus the score depends on the report only through the probability assigned to the
realized outcome $x$.
\end{definition}

Locality is economically natural: evaluation depends only on the probability
assigned to the event that occurs, not on probabilities assigned to counterfactual
outcomes. Log-likelihood evaluation is the canonical example.

\begin{definition}[Mutual information]
\label{def:mi}
Let $f(dy\mid s)$ be a channel with induced marginal $q(dy)$. Assume
$f(\cdot\mid s)\ll q$ for $P_S$-a.e.\ $s$. The mutual information between $S$ and $Y$ is
\begin{equation}
\label{eq:mi}
\MI(S;Y)
= \E_{P_{SY}}\!\left[\log\frac{df(\cdot\mid S)}{dq}(Y)\right]
= \int \KL\!\big(f(\cdot\mid s)\,\|\,q\big)\,P_S(ds).
\end{equation}
\end{definition}

The quantity $\MI(S;Y)$ measures how much reported predictions vary with the
information state. Larger mutual information corresponds to sharper state dependence
and therefore higher \emph{information complexity} of the rule.

\begin{definition}[Entropy-regularized objective]
\label{def:objective}
Fix a prior $P_S$, a utility function $U:\mathcal{S}\times\mathcal{Y}\to\R$, and a
regularization parameter $\lambda>0$. The entropy-regularized objective is
\begin{equation}
\label{eq:objective}
J_\lambda(f)
= \E_{P_{SY}}\!\big[U(S,Y)\big] - \lambda^{-1}\MI(S;Y).
\end{equation}
\end{definition}

The objective \eqref{eq:objective} trades off expected utility against information
complexity. Larger $\lambda$ permits more state-dependent rules; smaller $\lambda$
pushes toward state-insensitive reporting. The parameter $\lambda^{-1}$ is the
Lagrange multiplier associated with an information constraint $\MI(S;Y)\le \kappa$.

\section{Coherence and Uniqueness of the Logarithmic Score}
\label{sec:coherence}

We now show that the logarithmic score is not merely convenient, but is
\emph{required} by minimal coherence requirements on predictive evaluation.
The key requirement is \emph{aggregation coherence}: evaluation should be
invariant to whether one predicts at full resolution or predicts in stages
(coarse first, then fine).

\subsection{Coarsening, refinement, and additive evaluation}

Let $\{A_1,\ldots,A_K\}$ be a measurable partition of the outcome space $\mathcal{X}$.
For a predictive report $y\in\mathcal{Y}$, define the induced \emph{coarse} report
$y_{\mathrm{coarse}}(k) := y(A_k)$ and, whenever $y(A_k)>0$, the \emph{within-block}
conditional report $y_{\mathrm{fine}}(B\mid A_k) := y(B\cap A_k)/y(A_k)$.

\begin{assumption}[Aggregation coherence]
\label{ass:aggregation}
For any measurable partition $\{A_1,\ldots,A_K\}$, any predictive report $y$, and
any realized outcome $x\in A_k$ with $y(A_k)>0$, the score satisfies
\begin{equation}
\label{eq:aggregation}
S(x,y)
= S\!\big(k, y_{\mathrm{coarse}}\big)
+ S\!\big(x\mid A_k, y_{\mathrm{fine}}(\cdot\mid A_k)\big).
\end{equation}
\end{assumption}

Assumption~\ref{ass:aggregation} is the predictive analogue of Shannon's branching
axiom: evaluating a fine-resolution forecast equals evaluating the induced coarse
forecast plus the refinement forecast within the realized coarse cell.

Aggregation coherence has direct relevance for econometric practice. Many data
settings involve outcomes reported at varying levels of granularity: top-coded
income, grouped survey responses, interval-censored data, and administrative
records with binned variables. If evaluation rules are not coherent under such
coarsening, predictive performance comparisons become artifacts of reporting
conventions rather than substantive measures of forecast quality.

\subsection{Concrete interpretation of the coherence axioms}

The three axioms correspond to concrete requirements that applied econometricians
implicitly invoke when evaluating predictive performance.

\emph{Locality: ``Score what you predicted for what happened.''}
The score depends only on the probability assigned to the realized outcome.
If the model predicted $\hat{p}(Y_i = y_i \mid X_i) = 0.73$ and outcome $y_i$
materialized, the score depends on $0.73$, not on probabilities assigned to
alternative outcomes.
\emph{Strict Propriety: ``Honesty is uniquely optimal.''}
A forecaster maximizes expected score by reporting their true beliefs. Without
this, evaluation rules create incentives to shade probabilities strategically.
\emph{Aggregation Coherence: ``Resolution shouldn't matter.''}
Predicting ``which continent'' then ``which country within that continent''
should yield the same total score as predicting ``which country'' directly.

\begin{example}[Failure of aggregation coherence]
\label{ex:brier_fail}
The Brier (quadratic) score is strictly proper but fails aggregation coherence.
Consider outcomes $\{1,2,3\}$ with partition $\{\{1\}, \{2,3\}\}$ and true
distribution $(0.5, 0.3, 0.2)$. If outcome $2$ occurs, the direct Brier score
and the sum of coarse plus fine Brier scores do not match. The log score, by
contrast, satisfies $\log(0.3) = \log(0.5) + \log(0.6)$ exactly.
\end{example}

\subsection{Uniqueness of the log score under coherence}

Aggregation coherence distinguishes our characterization from Bregman families
of proper local scores, which generate divergences via convex generators without
refinement additivity.

\begin{theorem}[Uniqueness of the logarithmic score]
\label{thm:logscore}
Let $S:\mathcal{X}\times\mathcal{Y}\to\R$ satisfy (i) locality, (ii) strict
propriety, and (iii) aggregation coherence. Then there exist constants $c>0$
and a measurable function $d:\mathcal{X}\to\R$ such that
\begin{equation}
\label{eq:logscore}
S(x,y)=c\log\!\left(\frac{dy}{d\mu}(x)\right)+d(x).
\end{equation}
\end{theorem}

\paragraph{Proof sketch (full proof in \ref{app:logscore})}
Strict propriety together with locality forces a $1/p$-type derivative form for
the scoring function. Aggregation coherence imposes a multiplicative functional
equation corresponding to additivity under coarse--fine decomposition. Under mild
regularity, the unique solution is logarithmic. The result applies equally to
discrete, continuous, and mixed outcome spaces under standard measurability conditions.

\begin{corollary}[Mutual information as the coherent refinement gain]
\label{cor:kl}
Let $S$ be the logarithmic score and let $f(dy\mid s)$ be a stochastic predictive
rule with induced marginal $q(dy)$. The expected gain from state-dependent reporting
equals mutual information:
$\E\!\left[\log\frac{df(\cdot\mid S)}{dq}(Y)\right] = \MI(S;Y)$.
Equivalently, KL divergence is the unique complexity measure induced by coherent
local evaluation.
\end{corollary}

\noindent\emph{Takeaway.} Corollary~\ref{cor:kl} identifies KL/mutual information as the
unique refinement cost compatible with local, strictly proper, aggregation-coherent evaluation.
With the complexity term fixed, the next step is to solve the resulting regularized
optimization problem over predictive rules; \Cref{sec:gibbs} shows that the solution
necessarily takes a Gibbs/exponential-tilting form.

\section{Entropy-Regularized Predictive Inference and Gibbs Solutions}
\label{sec:gibbs}

We characterize the optimal predictive rule under the entropy-regularized objective.
The solution takes a Gibbs (exponential-tilt) form: this structure is not assumed
but entailed by the coherence-derived penalty.

\subsection{Gibbs form and self-consistency}

\begin{theorem}[Gibbs optimal rule and fixed point]
\label{thm:gibbs}
Under regularity conditions (\ref{app:regularity}), any maximizer $f^*$ of
$J_\lambda(f)$ admits the Gibbs--Boltzmann form
\begin{equation}
\label{eq:gibbs}
f^*(dy\mid s)
= \frac{\exp\{\lambda U(s,y)\}\, q^*(dy)}{Z_\lambda(s)},
\qquad
Z_\lambda(s)=\int_{\mathcal{Y}} \exp\{\lambda U(s,y')\}\,q^*(dy'),
\end{equation}
where $q^*$ is the marginal distribution of reports induced by $f^*$:
$q^*(dy)=\int_{\mathcal{S}} f^*(dy\mid s)\,P_S(ds)$.
The optimal value satisfies $J_\lambda(f^*)=\lambda^{-1}\E[\log Z_\lambda(S)]$.
If $U(s,\cdot)$ is not $q$-a.s.\ constant, the maximizer is unique up to $P_S$-null sets.
\end{theorem}

\paragraph{Proof sketch (full proof in \ref{app:gibbs})}
For fixed baseline $q$, variational calculus for relative-entropy regularization
yields the Gibbs form as the unique conditional maximizer. Self-consistency requires
$q^*$ to equal the induced marginal, establishing the fixed-point structure.

The Gibbs structure is transparent: for each state $s$, the optimal conditional
distribution over reports is an exponential tilt of the unconditional report
distribution $q^*$, with tilt magnitude controlled by $\lambda$. Larger $\lambda$
permits sharper state-dependent tilts; smaller $\lambda$ forces $f^*(\cdot\mid s)$
to remain close to $q^*$ in KL divergence.

\begin{remark}[Limiting behavior]
As $\lambda\to\infty$, the conditional law concentrates on utility-maximizing reports.
As $\lambda\to 0$, $f^*(dy\mid s)\Rightarrow q^*(dy)$ and $\MI(S;Y)\to 0$.
\end{remark}

\begin{remark}[Logit as entropy-regularized prediction]
For discrete choice among $K$ alternatives with utilities $(u_1,\ldots,u_K)$,
\eqref{eq:gibbs} yields $f^*(k\mid u)\propto \exp(\lambda u_k)$, i.e.\ multinomial
logit with inverse temperature $\lambda$. The logit form arises as a Gibbs solution,
not from assuming Gumbel shocks. This distinction matters empirically when baseline
probabilities are misspecified: the Gibbs tilt adjusts relative to the baseline,
and validation-selected $\eta$ adapts to baseline quality (see \ref{app:simulations}).
\end{remark}

\begin{remark}[Connection to IIA, Luce's axiom, and rational inattention]
\label{rem:luce}
Aggregation coherence is the predictive analogue of Luce's choice axiom
\citep{luce1959individual}: both impose separability under partitioning. The
logit/Gibbs form satisfies IIA precisely because the log score satisfies
aggregation coherence \citep{mcfadden1974conditional,yellott1977relationship}.
\citet{matejka2015rational} derive logit from a rational inattention problem
that \emph{assumes} Shannon mutual information as the cost. Our framework
\emph{derives} this cost from coherence: locality, strict propriety, and
aggregation coherence jointly force KL as the unique information measure,
making the logit structure a consequence of coherent evaluation rather than
an assumption.
\end{remark}

\subsection{Relation to Blahut--Arimoto}

It is useful to contrast Theorem~\ref{thm:gibbs} with the classical Blahut--Arimoto
algorithm for computing channel capacity and rate--distortion functions
\citep{blahut1972computation,arimoto1972algorithm}. The difference is ontological:
here we \emph{maximize} expected utility minus information-complexity, yielding a
Gibbs tilt toward high-utility reports, not a distortion-minimizing quantizer. The
fixed point $q^*$ is the unconditional distribution of predictive reports, not a
coding equilibrium. The result is not algorithmic but axiomatic: once KL/entropy
is selected by the coherence axioms, the Gibbs form follows.

\section{Complete-Class Characterization}
\label{sec:complete}

Having established that coherent evaluation uniquely selects KL divergence
(Section~\ref{sec:coherence}) and that optimal rules take the Gibbs form
(Section~\ref{sec:gibbs}), we now formalize the central structural claim:
\emph{Gibbs rules are the only admissible predictive rules} under a natural
ordering that values both predictive accuracy and informational parsimony.

\subsection{The efficiency frontier}

Consider a decision maker who cares about two objectives: expected predictive
utility $\E_f[U(S,Y)]$ and information complexity $\MI_f(S;Y)$. Higher utility
is better; lower complexity is better. A predictive rule $f$ \emph{dominates}
rule $g$ if $f$ achieves at least as much utility with at most as much complexity,
with strict inequality in at least one dimension. A rule is \emph{admissible}
if no other rule dominates it.

The set of admissible rules traces out an \emph{efficiency frontier} in
utility--complexity space. Points on this frontier represent optimal tradeoffs:
to gain more predictive accuracy, one must accept higher information costs.
Points below the frontier are wasteful: there exists a dominating rule that
predicts better or uses less information (or both).

\subsection{Gibbs rules span the frontier}

The key result is that the Gibbs family $\{f^*_\lambda : \lambda > 0\}$
parameterizes the entire efficiency frontier.

\begin{theorem}[Complete-class theorem]
\label{thm:complete}
Let $\mathcal{G} = \{f^*_\lambda : \lambda > 0\}$ denote the family of Gibbs
rules from Theorem~\ref{thm:gibbs}. Then:
\begin{enumerate}[label=(\roman*),nosep]
\item \emph{Gibbs rules are admissible.} Every $f^*_\lambda$ lies on the
efficiency frontier.

\item \emph{Non-Gibbs rules are dominated.} Any rule not of the Gibbs form
is strictly dominated by some Gibbs rule.

\item \emph{Essential completeness.} The Gibbs family is essentially complete
(in the sense of weak convergence on finite outcome spaces): every admissible
rule is either Gibbs or a mixture of two Gibbs rules at a frontier kink.
\end{enumerate}
\end{theorem}

\paragraph{Proof intuition (full proof in \ref{app:complete})}
Part (i): If $f^*_\lambda$ were dominated, some rule $g$ would achieve higher
$J_\lambda(g) = \E[U] - \lambda^{-1}\MI$, contradicting optimality of $f^*_\lambda$.
Part (ii): The frontier is the upper envelope of lines with slope $-\lambda^{-1}$.
Any point below this envelope is dominated by the frontier point directly above it.
Part (iii): At smooth frontier points, a unique Gibbs rule is tangent. At kinks,
admissible rules are convex combinations of adjacent Gibbs rules.

\emph{Scope.} The result applies to predictive rules mapping signals into
distributions. It does not directly address estimation strategies based on moment
conditions or estimating equations; these require additional structure linking
moments to predictive objectives, which we do not pursue here.

\subsection{Rational inattention as the boundary case}

The complete-class result clarifies the relationship between entropy-regularized
prediction and rational inattention.

\begin{corollary}[Rational inattention as constrained optimum]
\label{cor:ri}
Any solution to $\max_f \E_f[U(S,Y)]$ subject to $\MI_f(S;Y) \le \kappa$
coincides with a Gibbs rule $f^*_\lambda$, where $\lambda^{-1}$ is the shadow
price of the information constraint.
\end{corollary}

Rational inattention is not a separate modeling primitive but the boundary
regime of entropy-regularized prediction: it corresponds to points on the
efficiency frontier where the information budget binds.

\subsection{Practical implications}

The theorem has direct implications for applied work:

\emph{Use Gibbs rules.} Ad hoc predictive rules, ``shrink toward the prior
mean by 30\%'' or ``average MLE and OLS'', are typically inadmissible within the class we study.
There exists a Gibbs rule that dominates them.

\emph{The only choice is $\lambda$.} Once the Gibbs form is accepted, the
practitioner's task reduces to selecting the tradeoff parameter. Cross-validation
on a proper score (Section~\ref{app:practical}) provides a coherence-based selection
mechanism.

\emph{Coherence selects the functional form.} Unlike ridge or lasso, where
the penalty geometry is chosen exogenously, entropy regularization derives
its geometry from coherent evaluation. The complete-class theorem shows this
geometry is not merely convenient but uniquely characterized under the stated axioms.

\section{Identification Geometry and Local Asymptotic Interpretation}
\label{sec:identification}

This section formalizes the geometric mechanism by which entropy regularization
stabilizes the predictive criterion in weakly identified settings. The main message is that
the entropy penalty contributes curvature precisely in directions where the
unregularized criterion is flat, stabilizing the predictive map.

\subsection{The curvature problem in weak identification}

Weak identification corresponds to near-null directions in the criterion linking
parameters to observables: the Hessian (or Fisher information analogue) is singular
or ill-conditioned, Wald-type quadratic approximations become unreliable
\citep{moreira2003clr}, and small data perturbations induce large changes in
the optimizer.

The entropy penalty adds curvature in these directions. By contributing a term proportional to
mutual information, the regularized criterion acquires additional curvature
$\lambda^{-1}H_{\MI}(\theta)$. This contribution is positive semidefinite and,
under mild conditions, positive definite; it has mass in \emph{all} directions,
including those where the utility curvature $H_U(\theta)$ is deficient.

\subsection{Parametric restriction and regularized criterion}

Let $\{f_\theta:\theta\in\Theta\}$ be a parametric family of predictive rules with
$\Theta\subseteq\R^d$ open. The entropy-regularized objective restricted to this
family is
\begin{equation}
\label{eq:Jtheta}
J_\lambda(\theta)
= \E\!\left[U\!\bigl(X,Y;\theta\bigr)\right] - \lambda^{-1}\MI_\theta,
\end{equation}
with (negative) Hessian $H_\lambda(\theta)=-\nabla_\theta^2J_\lambda(\theta)$.
Identification is curvature of $J_\lambda$ near the optimizer; weak identification
corresponds to ill-conditioning of $H_\lambda(\theta^*_\lambda)$.

\subsection{Curvature decomposition and identification restoration}

\begin{proposition}[Curvature decomposition]
\label{prop:curvature}
Suppose $J_\lambda(\theta)$ is twice continuously differentiable. Then
\begin{equation}
\label{eq:hessian}
H_\lambda(\theta)
= H_U(\theta) + \lambda^{-1}H_{\MI}(\theta),
\end{equation}
where $H_U(\theta)=-\nabla_\theta^2\E[U(X,Y;\theta)]$ and
$H_{\MI}(\theta)=\nabla_\theta^2\MI_\theta$.
\end{proposition}

The decomposition shows that entropy regularization restores curvature whenever
$H_{\MI}(\theta)$ has mass in the directions where $H_U(\theta)$ is deficient.
This mechanism is structural: it does not depend on choosing an $L_1$ or $L_2$
penalty, but follows from the coherence-derived information penalty.

\emph{Weak-IV illustration.}
In the IV regression $Y_i = X_i\beta + \varepsilon_i$, $X_i = Z_i'\Pi + V_i$, weak
instruments ($\Pi\approx 0$) mean the curvature $H_U(\beta)$ in the $\beta$-direction
is small. The entropy penalty contributes $H_{\MI}(\beta)$, which measures how sharply
the predictive distribution responds to changes in $\beta$. Because this term is
positive whenever the predictive rule is nontrivial, it restores well-conditioning.

\subsection{Local asymptotic quadratic (LAQ) expansion}

We provide a local quadratic expansion as an interpretive bridge to classical
weak-identification diagnostics. The expansion is for the \emph{regularized
criterion} $J_{\lambda,n}$, not a log-likelihood ratio process; it is
criterion-level LAQ, not an automatic claim about the likelihood experiment.

\begin{assumption}[Local regularity]
\label{ass:lan}
At a pseudo-true maximizer $\theta_\lambda$, the utility is twice continuously
differentiable with finite variance of the gradient and bounded Hessian, and the
first-order condition $\nabla_\theta J_\lambda(\theta_\lambda)=0$ holds.
\end{assumption}

\begin{theorem}[LAQ expansion of the regularized criterion]
\label{thm:lan}
Under Assumption~\ref{ass:lan}, for local sequences $\theta_n=\theta_\lambda+h/\sqrt{n}$,
\begin{equation}
\label{eq:lan}
n\!\left(J_{\lambda,n}(\theta_n)-J_{\lambda,n}(\theta_\lambda)\right)
= h^\top \Delta_n
-\frac{1}{2}h^\top H_\lambda(\theta_\lambda) h + o_P(1),
\end{equation}
where $\Delta_n \Rightarrow N(0,\Sigma_\lambda)$ with
$\Sigma_\lambda = \Var(\nabla_\theta U(X,Y;\theta_\lambda))$ and
$H_\lambda(\theta_\lambda) = H_U(\theta_\lambda)+\lambda^{-1}H_{\MI}(\theta_\lambda)$.
\end{theorem}

\emph{Interpretation and scope.}
The expansion controls local geometry of the optimizer and supports curvature
diagnostics (eigenstructure, condition numbers, identification of near-flat
directions) without asserting LAN for the data-generating family. The key
econometric payoff is that when $H_\lambda$ is nonsingular, the expansion yields
standard local geometry for the regularized optimizer even when the unregularized
criterion exhibits weak identification: \emph{entropy contributes curvature in
the near-flat directions, stabilizing the predictive map}. The pseudo-true value
$\theta_\lambda$ maximizes the regularized predictive criterion, not a likelihood,
and coincides with the true parameter only under correct specification.

In the special case where $U$ is the log score and $\{f_\theta\}$ is a likelihood
family, the expansion aligns with LAN and the Hessian corresponds to regularized
Fisher information. Outside that setting, the LAQ statement remains a property of
the regularized criterion, with justification resting on coherence and strict
propriety rather than likelihood-based asymptotic normality.

\paragraph{Two research fronts and the bridge between them.}
The LAQ expansion identifies a connection between two largely independent
programs in weak-identification asymptotics.

On one front, \citet{andrewsmikusheva2016,andrewsmikusheva2022}, building on
\citet{mikusheva2007uniform,mikusheva2010robust}, develop a geometric approach
rooted in tangent cones and conditional inference. Weak identification is
characterized by the tangent space at the boundary of the identified set;
valid inference is achieved by conditioning on sufficient statistics for
nuisance parameters in the flat directions. This delivers LAMN limit
experiments and optimal decision rules when the degeneracy structure is known.

On a second front, \citet{antoine2009efficient,antoine2012efficient} develop
GMM-based asymptotics under nearly-weak identification where moment conditions
converge at rates between $n^{-1/2}$ and $n^{-1/4}$, and
\citet{antoine2020testing} propose identification-strength tests calibrated to
these intermediate regimes.

Our curvature decomposition $H_\lambda = H_U + \lambda^{-1}H_{\MI}$ connects
both fronts through predictive regularization. The Andrews--Mikusheva tangent
cone corresponds to $\ker(H_U)$: the flat directions of the unregularized
criterion. The entropy term $\lambda^{-1}H_{\MI}$ adds curvature precisely
there, and the regularized criterion remains well-conditioned across the
identification regimes that \citet{antoine2020testing} diagnostics would
distinguish. Entropy regularization thus smooths the discontinuous transitions
between regimes: rather than requiring different asymptotic approximations,
the regularized criterion admits a single LAQ structure throughout.

However, the LAQ expansion is deliberately limited in scope: it characterizes
local geometry but does not deliver distributional limits under drifting
identification. A complete asymptotic treatment, establishing convergence of
the regularized score process when the information matrix degenerates along
sequences, requires nonstandard probabilistic machinery (convergence of
degenerate diffusions, spectral methods for singular information operators)
that lies beyond the present paper. We identify the bridge but defer the
crossing; the curvature decomposition shows \emph{where} entropy regularization
acts, while the theory for \emph{how} the estimator behaves under drifting
identification necessitates unifying the geometric and GMM-based approaches.

\section{Monte Carlo Evidence: Weakly Identified IV}
\label{sec:numerical}

This section provides Monte Carlo evidence for the curvature mechanism
implied by entropy-regularized prediction, using the weak-instrument IV
setting introduced in Section~\ref{sec:intro}. We examine three
dimensions: how the mechanism varies with identification strength, its
robustness to baseline misspecification, and a visual decomposition of
curvature restoration. A supplementary discrete-choice simulation appears
in \ref{app:simulations}.

\subsection{Design}

We generate data from $Y_i = X_i\beta + \varepsilon_i$, $X_i = Z_i'\Pi + V_i$
with $Z_i \in \R^5$, $(\varepsilon_i, V_i)' \sim N(0,\Sigma)$, and
$\mathrm{Corr}(\varepsilon_i, V_i) = 0.5$. Instrument strength is
controlled by scaling $\Pi$ to achieve first-stage $F$-statistics of
approximately 3, 6, and 11, spanning severe, moderate, and
borderline-adequate identification. The sample size is $n = 400$
(50/25/25 training/validation/test split), with 200 Monte Carlo
replications. The entropy-regularized (E-R) predictor takes the Gibbs
form of Theorem~\ref{thm:gibbs}, with the information price
$\eta = \lambda^{-1}$ selected by validation on the held-out set.

\subsection{Instrument strength variation}

Table~\ref{tab:iv_strength} reports full results across three
identification regimes with a correctly specified baseline
($\beta_0 = \beta = 1$). The pattern is consistent with the curvature
decomposition in Proposition~\ref{prop:curvature}: as identification
weakens (lower $F$), the unregularized utility curvature $H_U(\beta)$
shrinks and both 2SLS and LIML deteriorate, the former through bias, the
latter through variance inflation. E-R stabilizes predictions across all
regimes by contributing entropy curvature
$\lambda^{-1}H_{\MI}(\beta)$ in the deficient direction.

\begin{table}[htbp]
\centering
\caption{Instrument strength variation ($\beta_0 = \beta = 1$, Normal errors).}
\label{tab:iv_strength}
\begin{tabular}{llccccccc}
\toprule
$F$ & Method & MSE & (SE) & Bias & SD & RMSE & Cov & $\bar{\eta}$ \\
\midrule
\multirow{3}{*}{$\approx 3$}
& 2SLS & 0.935 & (0.023) & 0.239 & 0.337 & 0.413 & 0.83 & ,  \\
& LIML & 3.821 & (1.014) & $-0.074$ & 1.630 & 1.632 & 0.84 & ,  \\
& E-R  & \emph{0.854} & (0.011) & 0.237 & 0.186 & 0.301 & 0.88 & 6.16 \\
\addlinespace
\multirow{3}{*}{$\approx 6$}
& 2SLS & 0.973 & (0.023) & 0.116 & 0.263 & 0.287 & 0.89 & ,  \\
& LIML & 1.656 & (0.259) & 0.032 & 0.767 & 0.768 & 0.90 & ,  \\
& E-R  & \emph{0.894} & (0.011) & 0.157 & 0.157 & 0.222 & 0.89 & 6.24 \\
\addlinespace
\multirow{3}{*}{$\approx 11$}
& 2SLS & 0.988 & (0.018) & 0.063 & 0.198 & 0.208 & 0.91 & ,  \\
& LIML & 1.083 & (0.025) & $-0.014$ & 0.234 & 0.235 & 0.93 & ,  \\
& E-R  & \emph{0.925} & (0.011) & 0.106 & 0.125 & 0.164 & 0.91 & 6.37 \\
\bottomrule
\end{tabular}
\begin{flushleft}
\small\textit{Notes:} MSE is out-of-sample predictive mean squared error;
Cov is 95\% CI coverage; $\bar\eta$ is validation-selected information
price. $n = 400$, 200 replications.
\end{flushleft}
\end{table}

Under severe weakness ($F \approx 3$), E-R reduces predictive MSE by
approximately 9\% relative to 2SLS (0.854 vs.\ 0.935). The improvement
is driven by variance reduction (SD of 0.186 vs.\ 0.337) rather than bias
correction, the signature of curvature restoration. LIML, which lacks
regularization, exhibits a standard deviation of 1.63, roughly five times
that of 2SLS, making it unreliable for prediction despite its asymptotic
unbiasedness. As identification strengthens ($F \approx 11$), the
three methods converge, consistent with the LAQ expansion
(Theorem~\ref{thm:lan}): when $H_U$ is already well-conditioned, the
entropy contribution is proportionally modest.

\subsection{Robustness to baseline misspecification}

A natural concern is sensitivity to the baseline $\beta_0$. We repeat
the experiment with $\beta_0 = 0$ (misspecified; the true $\beta = 1$).

\begin{table}[htbp]
\centering
\caption{Misspecified baseline ($\beta_0 = 0$, $\beta = 1$, Normal errors).}
\label{tab:iv_misspec}
\begin{tabular}{llccccccc}
\toprule
$F$ & Method & MSE & (SE) & Bias & SD & RMSE & Cov & $\bar{\eta}$ \\
\midrule
\multirow{3}{*}{$\approx 3$}
& 2SLS & 0.935 & (0.023) & 0.239 & 0.337 & 0.413 & 0.83 & ,  \\
& LIML & 3.821 & (1.014) & $-0.074$ & 1.630 & 1.632 & 0.84 & ,  \\
& E-R  & \emph{0.931} & (0.024) & 0.186 & 0.280 & 0.336 & 0.90 & 8.50 \\
\addlinespace
\multirow{3}{*}{$\approx 6$}
& 2SLS & 0.973 & (0.023) & 0.116 & 0.263 & 0.287 & 0.89 & ,  \\
& LIML & 1.656 & (0.259) & 0.032 & 0.767 & 0.768 & 0.90 & ,  \\
& E-R  & 0.978 & (0.023) & 0.100 & 0.249 & 0.268 & 0.90 & 9.25 \\
\bottomrule
\end{tabular}
\begin{flushleft}
\small\textit{Notes:} Baseline centered at $\beta_0 = 0$; true parameter
$\beta = 1$. Other settings as in Table~\ref{tab:iv_strength}.
\end{flushleft}
\end{table}

Table~\ref{tab:iv_misspec} shows that validation-based selection of $\eta$
adapts to misspecification: the selected $\bar\eta$ increases (8.50 vs.\
6.16 under $F \approx 3$), reflecting the optimizer's reluctance to shrink
toward an inaccurate baseline. Under severe weakness, E-R still matches
2SLS (0.931 vs.\ 0.935) without the catastrophic variance of LIML. Under
moderate weakness ($F \approx 6$), E-R and 2SLS are essentially tied.
The mechanism is transparent: validation selects a larger $\eta$ when the
baseline is wrong, tilting less aggressively and preserving the
data-driven signal. This self-correcting behavior is a direct consequence
of the Gibbs structure (Theorem~\ref{thm:gibbs}): the exponential tilt
adjusts relative to the baseline, and validation selects the tilt
magnitude that optimizes out-of-sample prediction.

\subsection{Curvature restoration: visual evidence}

\begin{figure}[htbp]
\centering
\includegraphics[width=\textwidth]{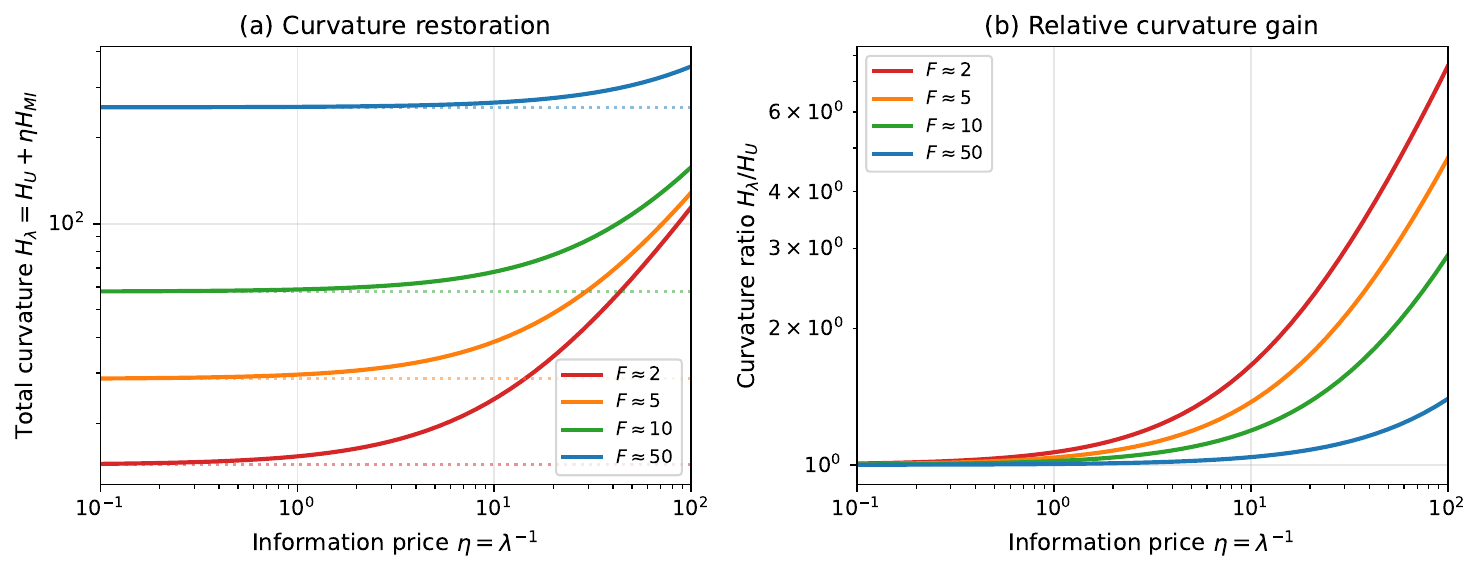}
\caption{Curvature restoration under entropy regularization. (a) Total curvature
$H_\lambda = H_U + \eta H_{\MI}$ as a function of information price $\eta$ for
different first-stage $F$-statistics; dotted lines show unregularized curvature $H_U$.
(b) Curvature ratio $H_\lambda / H_U$ showing relative gain from entropy penalty.
Under weak instruments ($F \approx 2$), entropy contributes most of the curvature;
under strong instruments ($F \approx 50$), the data-driven curvature dominates.
$n = 400$, $k = 5$, $\rho = 0.5$, averaged over 50 replications.}
\label{fig:curvature}
\end{figure}

Figure~\ref{fig:curvature} provides a direct visualization of the curvature
decomposition $H_\lambda = H_U + \eta H_{\MI}$ from
Proposition~\ref{prop:curvature}.

Panel~(a) shows total curvature as a function of $\eta$ for different
first-stage $F$-statistics. Under weak instruments ($F \approx 2$), the
unregularized curvature $H_U$ (dotted line) is very low; as $\eta$
increases, the entropy term $\eta H_{\MI}$ adds substantial curvature.
Under strong instruments ($F \approx 50$), $H_U$ is already large and
the relative contribution is modest. Panel~(b) quantifies this: for
$F \approx 2$, curvature is amplified roughly tenfold at $\eta = 100$,
while for $F \approx 50$ the amplification is less than twofold.
This visualizes the core geometric mechanism: entropy regularization
matters most precisely when identification is weakest.

\section{Discussion}
\label{sec:discussion}

This paper contributes a perspective in econometrics from estimator-centric analysis
to \emph{prediction-centric} evaluation under limited identification. The answer
to ``what penalty?'' is not an algorithmic choice but a structural one: once
prediction is evaluated coherently, entropy regularization and the Gibbs form
are structurally implied.

\subsection{Entropy regularization, Bregman geometry, and what is special about KL}

Modern econometric practice embraces a wide class of regularizers, $\ell_1$,
$\ell_2$, elastic net, typically justified by sparsity, tractability, or oracle
bounds, imposing geometry in \emph{parameter space}. Entropy regularization differs:
it imposes geometry in \emph{predictive space}, measured by divergence between
predictive distributions.

This distinction becomes sharp when comparing KL to Bregman or $f$-divergences.
Alternative information costs may capture behavioral flexibility or non-logit
choice patterns \citep{FosgerauEtAl2020,BloedelDentiPomatto2025}, but recent results
show such divergences fail \emph{global coherence}: they are not Blackwell monotone
and cannot serve as globally valid measures of predictive refinement. KL is special
precisely because it is local, additive, and globally coherent.

Bregman penalties can function as legitimate \emph{local} regularizers: around an
interior optimum, they supply stabilizing curvature, but they lack a global
interpretation as admissible restrictions on predictive refinement. They should be
viewed as local approximations to coherent entropy regularization.

\subsection{Weak identification, LAQ humility, and beyond oracle asymptotics}

Our analysis reframes weak identification geometrically: it corresponds to flat
predictive criteria, and instability arises because unrestricted refinement amplifies
noise in flat directions. The curvature decomposition
$H_\lambda(\theta) = H_U(\theta) + \eta H_{\MI}(\theta)$
shows that stabilization arises from the information term, not distributional
assumptions (Figure~\ref{fig:curvature}). The LAQ expansion is deliberately limited
in scope: it explains why prediction stabilizes without asserting classical asymptotic
normality. As discussed in Section~\ref{sec:identification}, the curvature
decomposition links two largely independent research fronts: the geometric/tangent-cone
approach of \citet{andrewsmikusheva2016,andrewsmikusheva2022} and the GMM-based
nearly-weak asymptotics of \citet{antoine2009efficient,antoine2020testing}. Entropy
regularization provides a predictive mechanism that smooths across identification
regimes, but the complete asymptotic theory for drifting identification requires
nonstandard machinery that we defer to future work.

\subsection{Partial identification and policy evaluation}

Entropy-regularized prediction offers a disciplined selection rule over partially
identified sets of predictive distributions. The information price $\eta$ plays a
role analogous to a policy maker's tolerance for ambiguity: larger $\eta$ yields
more conservative predictions, while smaller $\eta$ permits sharper counterfactual
statements when supported by the data. Unlike minimax criteria, entropy
regularization preserves probabilistic structure and remains compatible with proper
scoring and out-of-sample evaluation.

\subsection{Implications for applied econometrics}

The analysis yields several concrete implications for practice:

\begin{itemize}[nosep]
\item \emph{Regularize in predictive space, not parameter space.} Under weak
identification, prediction should be stabilized by penalizing informational
complexity of the predictive map, not by shrinking coefficients toward zero.

\item \emph{Within our axioms, KL/entropy is implied by coherence rather than chosen for convenience.}
The penalty is uniquely selected by evaluation axioms, not chosen for
computational convenience or robustness.

\item \emph{The information price $\eta$ is diagnostic, not just tuning.}
Validation-selected $\eta > 5$ signals weak identification; $\eta < 1$ signals
strong data support for sharp predictions.
\end{itemize}
Additional discussion on practical recommendations is presented in \ref{app:practical}.

\subsection{Concluding remarks}

The axiomatic starting point is deliberate. Rather than selecting a divergence
for computational convenience or asymptotic optimality and then exploring its
consequences, we ask what penalty on predictive refinement is \emph{implied} by
minimal requirements on how predictions are evaluated. The answer, entropy
regularization with KL complexity and the Gibbs form, is not a modeling choice but a
structural implication of locality, strict propriety, and aggregation coherence. The complete-class theorem confirms this within the admissible class: every non-Gibbs
rule is dominated.

The curvature decomposition $H_\lambda = H_U + \lambda^{-1}H_{\MI}$ translates
these axiomatics into a geometric mechanism: entropy contributes curvature in
flat directions, stabilizing prediction precisely where identification weakens.
This mechanism places in correspondence two largely independent programs in the weak-identification
literature. On the geometric side, the tangent-cone approach of
\citet{andrewsmikusheva2016,andrewsmikusheva2022} characterizes identification
failure through the structure of the information operator and delivers LAMN
limits via conditional inference. On the GMM side, the nearly-weak framework
of \citet{antoine2009efficient,antoine2020testing} allows moment conditions to
vanish at intermediate rates and develops tests calibrated to these regimes. Our
curvature decomposition shows that entropy regularization smooths across both:
the kernel of $H_U$ is the Andrews--Mikusheva tangent cone, and the entropy
term stabilizes the criterion across the identification regimes that
Antoine--Renault diagnostics distinguish. The LAQ expansion identifies this
bridge; a complete distributional theory under drifting identification, unifying
the geometric and moment-condition approaches, remains an important open problem
requiring nonstandard probabilistic machinery.

Econometric instabilities (weak identification, erratic finite-sample behavior,
overfitting) can be viewed as failures to control predictive refinement coherently.
By grounding regularization in the same coherence requirements that justify the
scoring rules used to evaluate it, the axiomatic approach ensures that the cure
is compatible with the diagnosis.
\bibliographystyle{apalike}
\bibliography{references}

\newpage
\appendix

\section{Technical Proofs and Regularity}
\label{app:proofs}

\subsection{Standing regularity conditions}
\label{app:regularity}

\begin{assumption}[Measurability and standard Borel structure]
\label{ass:measurability}
The spaces $(\mathcal{X},\mathcal{F}_X)$ and $(\mathcal{Y},\mathcal{F}_Y)$ are
standard Borel spaces. For each $x\in\mathcal{X}$, the predictive rule
$f(\cdot\mid x)$ is a probability measure on $(\mathcal{Y},\mathcal{F}_Y)$, and
$x\mapsto f(B\mid x)$ is measurable for every $B\in\mathcal{F}_Y$. The utility
function $U:\mathcal{X}\times\mathcal{Y}\to\mathbb{R}$ is jointly measurable.
\end{assumption}

\begin{assumption}[Integrability and finite partition functions]
\label{ass:integrability}
For the range of $\lambda$ under consideration,
$\E[\exp\{\lambda |U(X,Y)|\}] < \infty$ under the joint law induced by $(P_X,f)$.
\end{assumption}

\begin{assumption}[Interchange for parametric restrictions]
\label{ass:interchange}
When a parametric restriction $\{f_\theta\}$ is imposed, assume
$\theta\mapsto \E[U(X,Y;\theta)]$ and $\theta\mapsto \MI_\theta$ are twice
continuously differentiable, and differentiation may be interchanged with integration.
\end{assumption}

\subsection{Proof of Theorem~\ref{thm:logscore} (Uniqueness of the log score)}
\label{app:logscore}

We provide a discrete-space proof isolating two restrictions: strict propriety
under locality, and aggregation coherence.

Assume $\mathcal{X}=\{1,\ldots,K\}$ and a local scoring rule $S(i,p)=s(p_i)$ for
some function $s:(0,1]\to\mathbb{R}$.

\paragraph{Step 1 (Strict propriety + locality)}
Let $K=2$ with true distribution $p=(p,1-p)$ and report $q=(q,1-q)$. Propriety
requires $\phi_p(q)=p\,s(q)+(1-p)s(1-q)$ to be uniquely maximized at $q=p$.
Differentiating yields $p\,s'(p)=(1-p)\,s'(1-p)$ for all $p\in(0,1)$, hence
$s'(t)=c/t$ for some constant $c>0$ (strict propriety forces $c>0$). Integrating
yields $s(t)=c\log t + d$.

\paragraph{Step 2 (Aggregation coherence)}
Consider $\mathcal{X}=\{1,2,3\}$ and partition $\{\{1\},\{2,3\}\}$. Aggregation
coherence implies $s(rt)=s(r)+s(t)+C$. Substituting $s(t)=c\log t + d$ verifies
this equation; continuity plus Step 1 implies uniqueness.
\hfill$\square$

\subsection{Proof of Corollary~\ref{cor:kl}}
\label{app:kl}

Under the log score, $\E[S(X,Y)]=\E[\log f(Y\mid X)] + \text{const}$. The expected
improvement from using state-dependent $f(\cdot\mid X)$ rather than baseline $q$ is
$\E[\log f(Y\mid X)/q(Y)] = \MI(X;Y)$.
\hfill$\square$

\subsection{Proof of Theorem~\ref{thm:gibbs} (Gibbs optimal rule)}
\label{app:gibbs}

Write $\MI(X;Y)=\E[\log f(Y\mid X)]-\E[\log q(Y)]$. Fixing $x$ and optimizing over
$f(\cdot\mid x)$ subject to normalization yields the first-order condition
$\log f(y\mid x) = \lambda U(x,y)+\log q(y)-\log Z_\lambda(x)$, which is
\eqref{eq:gibbs}. Uniqueness follows from strict concavity.
\hfill$\square$

\subsection{Proof of Proposition~\ref{prop:curvature} (Curvature decomposition)}
\label{app:curvature}

Differentiate $J_\lambda(\theta)=\E[U(X,Y;\theta)]-\lambda^{-1}\MI_\theta$ twice:
$H_\lambda(\theta) = -\nabla_\theta^2\E[U] + \lambda^{-1}\nabla_\theta^2\MI_\theta
= H_U(\theta)+\lambda^{-1}H_{\MI}(\theta)$.
\hfill$\square$

\subsection{Proof of Theorem~\ref{thm:lan} (LAQ expansion)}
\label{app:lan}

Taylor expand the sample utility around $\theta_\lambda$. The CLT yields
$n^{-1/2}\sum_i \nabla_\theta U(X_i,Y_i;\theta_\lambda) \Rightarrow N(0,\Sigma_\lambda)$.
The law of large numbers handles the Hessian. The information term contributes
deterministic second-order curvature. Combining yields \eqref{eq:lan}.
\hfill$\square$

\section{Proofs for Complete-Class Theorem}
\label{app:complete}

This appendix provides formal definitions and complete proofs for the
complete-class characterization stated in Section~\ref{sec:complete}.

\subsection{Formal definitions}

\begin{definition}[Risk--complexity dominance]
\label{def:dominance}
For predictive rules $f$ and $g$, write $f \succ g$ if
$\E_f[U(S,Y)] \ge \E_g[U(S,Y)]$ and $\MI_f(S;Y) \le \MI_g(S;Y)$, with at least
one inequality strict. A rule is \emph{admissible} if it is not dominated.
\end{definition}

\begin{definition}[Value function and frontier]
Define the value function
$V(\kappa) = \sup\{\E_f[U(S,Y)] : \MI_f(S;Y) \le \kappa\}$.
The set $\{(\kappa, V(\kappa))\}$ is the efficient frontier.
\end{definition}

\begin{lemma}[Frontier properties]
\label{lem:frontier}
The value function $V(\kappa)$ is nondecreasing and concave on its effective domain.
\end{lemma}

\begin{proof}
Nondecreasing: If $\kappa' > \kappa$, the constraint set for $V(\kappa')$ contains
that for $V(\kappa)$. Concavity: For $\alpha \in (0,1)$ and rules $f$, $g$ achieving
$V(\kappa_1)$, $V(\kappa_2)$, the mixture $\alpha f + (1-\alpha)g$ is feasible for
$\alpha\kappa_1 + (1-\alpha)\kappa_2$ (by convexity of MI in the rule) and achieves
utility $\alpha V(\kappa_1) + (1-\alpha)V(\kappa_2)$.
\end{proof}

\subsection{Proof of Theorem~\ref{thm:complete}}

\paragraph{(i) Admissibility of Gibbs rules}
Suppose $f^*_\lambda$ is dominated by some rule $g$. Then either
$\E_g[U] > \E_{f^*_\lambda}[U]$ with $\MI_g \le \MI_{f^*_\lambda}$, or
$\E_g[U] \ge \E_{f^*_\lambda}[U]$ with $\MI_g < \MI_{f^*_\lambda}$.
In either case,
\[
J_\lambda(g) = \E_g[U] - \lambda^{-1}\MI_g
> \E_{f^*_\lambda}[U] - \lambda^{-1}\MI_{f^*_\lambda} = J_\lambda(f^*_\lambda),
\]
contradicting optimality of $f^*_\lambda$ for $J_\lambda$.

\paragraph{(ii) Gibbs characterization at smooth points}
At any frontier point $(\kappa, V(\kappa))$ where $V$ is differentiable,
there exists a unique supporting hyperplane with slope $V'(\kappa)$.
The Lagrangian for maximizing utility subject to $\MI \le \kappa$ yields
first-order conditions identical to those for $J_\lambda$ with
$\lambda^{-1} = V'(\kappa)$. By Theorem~\ref{thm:gibbs}, the solution
is the Gibbs rule $f^*_\lambda$.

\paragraph{(iii) Essential completeness}
At nondifferentiable points (frontier kinks), the subdifferential
$\partial V(\kappa)$ is an interval $[\lambda_1^{-1}, \lambda_2^{-1}]$.
Any admissible rule at such a point can be expressed as a convex combination
of $f^*_{\lambda_1}$ and $f^*_{\lambda_2}$. The topology for ``essentially
complete'' is weak convergence on $\Delta(\mathcal{Y})$; on finite outcome
spaces, this coincides with total variation convergence.
\hfill$\square$

\subsection{Proof of Corollary~\ref{cor:ri}}

The constrained problem $\max_f \E_f[U]$ s.t.\ $\MI_f \le \kappa$ has Lagrangian
$\mathcal{L} = \E[U] - \mu(\MI - \kappa)$. At an interior optimum, $\mu > 0$
and the first-order conditions match those for $J_\lambda$ with $\lambda^{-1} = \mu$.
The shadow price $\mu \in \partial V(\kappa)$ is the marginal value of relaxing
the information constraint.
\hfill$\square$

\section{Supplementary Simulations: Discrete Choice}
\label{app:simulations}

This appendix reports supplementary Monte Carlo results for discrete choice,
complementing the IV simulations in Section~\ref{sec:numerical}.

\subsection{Design and metrics}

For discrete choice we report negative log-likelihood (NLL), Brier score,
and expected calibration error (ECE). Monte Carlo standard errors appear
in parentheses.

\subsection{Discrete Choice with Sparse Partworths}

\emph{Design.}
Multinomial choice among $K=4$ alternatives with $U_{ik}=X_{ik}'\beta_k+\varepsilon_{ik}$
and Type-I extreme value errors. Parameters: $p=4$ attributes, 150/50/80 split, 50 replications.

\emph{E--R rule.}
Temperature scaling: $P^*_\eta(k\mid x) = \hat p_k(x)^{1/\eta}/\sum_{j=1}^K \hat p_j(x)^{1/\eta}$.

\begin{table}[htbp]
\centering
\caption{Discrete choice with sparse partworths: out-of-sample metrics.}
\begin{tabular}{llcccccc}
\toprule
Sparsity & Method & NLL & (SE) & Brier & (SE) & ECE & $\bar{\eta}$ \\
\midrule
\multirow{3}{*}{0\%}
& MLE   & 0.965 & (0.024) & 0.523 & (0.013) & 0.109 & ,  \\
& Ridge & \emph{0.955} & (0.023) & \emph{0.519} & (0.012) & \emph{0.105} & ,  \\
& E-R   & 0.964 & (0.023) & 0.522 & (0.012) & 0.108 & 1.17 \\
\addlinespace
\multirow{3}{*}{40\%}
& MLE   & 1.024 & (0.019) & 0.557 & (0.010) & 0.107 & ,  \\
& Ridge & \emph{1.017} & (0.018) & \emph{0.555} & (0.010) & \emph{0.110} & ,  \\
& E-R   & 1.031 & (0.019) & 0.560 & (0.010) & 0.113 & 1.18 \\
\addlinespace
\multirow{3}{*}{80\%}
& MLE   & 1.268 & (0.017) & 0.686 & (0.009) & 0.108 & ,  \\
& Ridge & \emph{1.260} & (0.015) & \emph{0.683} & (0.008) & \emph{0.090} & ,  \\
& E-R   & 1.262 & (0.014) & 0.684 & (0.008) & 0.096 & 1.56 \\
\bottomrule
\end{tabular}
\end{table}

\subsection{Discussion of discrete-choice results}

The results show that entropy-regularized prediction performs comparably to
ridge regularization across all sparsity levels, but does not offer substantial
finite-sample improvement. This is unsurprising and consistent with the broader
regularization literature in discrete choice. In well-specified logit models,
the likelihood itself possesses information-geometric structure (the logit form
is already a Gibbs solution; see Remark~\ref{rem:luce}), and parameter-space
regularization is effective because the model geometry is well-behaved.
Ridge regression, in particular, corresponds to a Gaussian prior on coefficients,
and efficient Bayesian computation for logistic models is available via
P\'{o}lya--Gamma data augmentation \citep{polson2013bayesian}, which provides
exact posterior sampling without Metropolis--Hastings tuning.

More broadly, the local-global shrinkage literature
\citep{polson2010shrink,carvalho2010horseshoe,polson2012local} demonstrates
that adaptive parameter-space regularization achieves strong predictive
performance in high-dimensional settings, including discrete choice, by
concentrating shrinkage on noise directions while preserving signal. These
methods exploit the specific structure of the parameter-space geometry and are
well-suited to settings where identification is uniform across directions.

The contribution of entropy regularization is therefore not finite-sample
predictive superiority in well-specified discrete choice, but rather the
coherence foundation: it is the unique penalty compatible with aggregation-coherent
evaluation (Theorem~\ref{thm:logscore}), and its advantages become decisive
precisely when identification degenerates directionally, as in the IV setting
of Section~\ref{sec:numerical}, where parameter-space penalties lack the geometric
structure to target the deficient directions.

\section{Practical Guide}
\label{app:practical}

This appendix provides practical guidance for implementing entropy-regularized
predictive inference. The goal is to translate the theoretical framework into
actionable recommendations for applied work.

\subsection{When to use entropy-regularized predictive inference}

Entropy regularization addresses a specific failure mode: instability arising from
unrestricted predictive refinement when the criterion linking parameters to
observables lacks curvature. The method is most valuable when identification is
weak or uncertain, manifesting as large standard errors, sensitivity to minor
specification changes, low first-stage $F$-statistics in IV, or nearly flat
likelihood surfaces. It is also well-suited to settings where the goal is
prediction rather than parameter recovery, where calibration matters (credit
scoring, medical diagnosis, demand forecasting), and where the baseline or prior
is informative but uncertain.

Conversely, entropy regularization is less useful when identification is strong
and samples are large, when the target is a point-identified structural parameter
with clear economic interpretation, or when coefficient sparsity is the primary
concern. In the latter case, $\ell_1$-type penalties may be more appropriate for
variable selection, though they lack the coherence properties of KL regularization.

\subsection{Choosing the information price \texorpdfstring{$\eta = \lambda^{-1}$}{eta = 1/lambda}}

The information price $\eta = \lambda^{-1}$ governs the tradeoff between predictive
sharpness and stability. Larger $\eta$ yields more conservative predictions that
deviate less from the baseline; smaller $\eta$ permits sharper state dependence
when supported by the data.

\emph{Recommended procedure: Proper-score cross-validation.}
The following procedure is robust across settings:
\begin{enumerate}[nosep]
\item Split data into training, validation, and test sets (e.g., 50/25/25).
\item For a grid of $\eta$ values (e.g., $\eta \in \{0.1, 0.5, 1, 2, 5, 10, 20\}$),
fit on training data and evaluate Gibbs-tilted predictions on validation using
a proper score (log score for density forecasts, Brier score for probabilities).
\item Select $\eta^*$ minimizing validation loss.
\item Report test-set performance at $\eta^*$.
\end{enumerate}

\emph{Interpretation of selected $\eta$.}
The validation-selected $\eta$ is itself diagnostic. If $\eta^*$ is small
(say, $< 1$), the data support sharp state dependence and regularization
contributes little. If $\eta^*$ is large (say, $> 5$), the data do not support
fine predictive refinement, and predictions should be conservative. Very large
$\eta^*$ signals potential weak identification: the curvature contribution
from the entropy term dominates the utility curvature.

\emph{Alternative: Marginal likelihood or information criteria.}
When cross-validation is computationally expensive, marginal likelihood (in the
Bayesian interpretation) or information criteria (AIC, BIC adapted for the
regularized criterion) provide alternatives. However, proper-score cross-validation
has the advantage of directly optimizing predictive performance on held-out data.

\subsection{Diagnosing weak identification}

Several diagnostics suggest that entropy regularization may help stabilize
predictions:
\begin{enumerate}[label=(\alph*),nosep]
\item First-stage $F < 10$ in IV regression (the Stock--Yogo threshold).
\item Condition number of the Hessian exceeds 100.
\item Large discrepancy between 2SLS and LIML estimates, or high LIML variance.
\item Validation-selected $\eta > 5$.
\item High bootstrap variance of predictions relative to parameter estimates.
\item Sensitivity of predictions to minor specification changes (adding/dropping
instruments, changing functional form).
\end{enumerate}

\emph{A note on LIML under weak identification.}
LIML removes asymptotic bias by using a larger correction factor ($\kappa \geq 1$),
but this comes at the cost of substantially increased variance when instruments
are weak. In our simulations with $F \approx 3$, LIML achieves near-zero bias
($-0.07$) but exhibits a standard deviation of 1.63, roughly five times that of
2SLS (0.34). This variance inflation makes LIML unreliable for prediction even
though it is asymptotically unbiased. The predictive MSE of LIML (3.82) far exceeds
that of 2SLS (0.94) and E-R (0.85). This illustrates the key insight: for prediction
under weak identification, bias-variance tradeoff matters more than asymptotic
unbiasedness. Entropy regularization achieves a favorable tradeoff by contributing
curvature that stabilizes predictions without the extreme variance of LIML.

When multiple diagnostics point toward weak identification, entropy regularization
provides a coherence-based approach to stabilize predictions without ad hoc shrinkage
toward zero or arbitrary prior centering.

\subsection{Software implementation}

For IV regression, the entropy-regularized predictor admits a closed-form
expression as a weighted average of the 2SLS estimate and a baseline:
\[
\hat\beta_{\eta} = \frac{(X'P_Z X/\sigma^2)\hat\beta_{\text{2SLS}}
+ (\eta/\tau^2)\beta_0}{X'P_Z X/\sigma^2 + \eta/\tau^2},
\]
where $P_Z = Z(Z'Z)^{-1}Z'$ is the instrument projection matrix, $\beta_0$ is the
baseline (prior mean), and $\tau^2$ is the baseline variance. This can be
implemented in any statistical software with matrix operations.

For discrete choice, the Gibbs-tilted predictor applies temperature scaling to
baseline probabilities:
\[
P^*_\eta(k\mid x) = \frac{\hat p_k(x)^{1/\eta}}{\sum_{j=1}^K \hat p_j(x)^{1/\eta}}.
\]
This requires only the baseline predicted probabilities and a single tuning
parameter. Code implementing these methods is available from the authors upon
request.

\end{document}